\documentclass[11pt]{article}
\usepackage{amsfonts}
\usepackage{mathrsfs}
\usepackage{bbm}
\usepackage{amsfonts}
\usepackage{amssymb,amsmath,graphicx}
\usepackage{float} \usepackage[colorlinks=true]{hyperref}
\hypersetup{urlcolor=blue, citecolor=red}

\allowdisplaybreaks

\begin{document}
\title{Global boundedness of solutions in a parabolic-parabolic chemotaxis
 system with singular sensitivity }
\author{Xiangdong Zhao \quad Sining Zheng{\thanks{Corresponding author.
E-mail: snzheng@dlut.edu.cn (S. N. Zheng), 704456001@qq.com (X. D. Zhao)}}\\
\footnotesize School of Mathematical Sciences, Dalian University of
Technology, Dalian 116024, P. R. China} \maketitle
\date{}
\newtheorem{theorem}{Theorem}
\newtheorem{definition}{Definition}[section]
\newtheorem{lemma}{Lemma}[section]
\newtheorem{proposition}{Proposition}[section]
\newtheorem{corollary}{Corollary}[section]
\newtheorem{remark}{Remark}
\renewcommand{\theequation}{\thesection.\arabic{equation}}
\catcode`@=11 \@addtoreset{equation}{section} \catcode`@=12
\maketitle{}
\begin{abstract}

We consider a parabolic-parabolic Keller-Segel system of chemotaxis
model with singular sensitivity $u_t=\Delta u-\chi\nabla\cdot(\frac{u}{v}\nabla
v)$, $v_t=k\Delta v-v+u$ under homogeneous Neumann boundary
conditions in a smooth bounded domain $\Omega\subset\mathbb{R}^n$
$(n\geq2)$, with $\chi,k>0$. It is proved that for any $k>0$, the
problem admits global classical solutions, whenever
$\chi\in\big(0,-\frac{k-1}{2}+\frac{1}{2}\sqrt{(k-1)^2+\frac{8k}{n}}\big)$.
The global solutions are moreover globally bounded if $n\le 8$. This
shows an exact way the size of the diffusion constant $k$ of the
chemicals $v$ effects the behavior of solutions.
\begin{description}
\item[2010MSC:] 35B40; 92C17; 35K55
\item[Keywords:] Chemotaxis; Logarithmic sensitivity; Boundedness
\end{description}
\end{abstract}


\section{Introduction}
In this paper, we consider the parabolic-parabolic chemotaxis system with singular sensitivity
\begin{eqnarray}\label{P} \left\{
 \begin{array}{llll}
  \displaystyle u_t=\Delta u-\chi\nabla \cdot(\frac{u}{v}\nabla v),& x\in\Omega,~~t\in(0,T),\\[4pt]
  \displaystyle v_t=k\Delta v-v+u,& x\in\Omega,~~t\in(0,T),\\[4pt]
  \displaystyle \frac{\partial u}{\partial {\nu}}=\frac{\partial v}{\partial {\nu}}=0 ,& x\in\partial\Omega,~~t\in(0,T),\\[4pt]
  \displaystyle (u(x,0),v(x,0))=(u_0(x),v_0(x)),  &x\in\Omega,
 \end{array}\right.
\end{eqnarray}
where $\chi,k>0$, $\Omega$ is a smooth bounded domain in
$\mathbb{R}^n~(n\geq 2)$, $\frac{\partial}{\partial \nu}$ denotes the derivation with respect to the outer normal of $\partial\Omega$, and the initial datum $u_0\in C^0(\overline{\Omega})$, $u_0(x)\geq 0$ on $\overline\Omega$, $v_0\in W^{1,q}(\Omega)$ $(q>n)$, $v_0(x)>0$ on $\overline\Omega$.

The classical Keller-Segel system of chemotaxis model was introduced by Keller
and Segel \cite{EL} in 1970 to describe the cells (with density $u$) move towards
the concentration gradient of a chemical substance $v$ produced by the cells themselves.
 Various forms of sensitivity functions can be chosen to model different types of
chemotaxis mechanisms. Among them $\phi(v)=\frac{\chi}{v}$ was
selected in \eqref{P} largely due to the Weber-Fechner's law for
cellular behaviors, where the subjective sensation is proportional
to the logarithm of the stimulus intensity \cite{TK}. With
$\phi(v)=\frac{\chi}{v}$, the cellular movements are governed by the
taxis flux $\frac{\chi\nabla v}{v}$, which may be unbounded when
$v\approx 0$.  In the model \eqref{P}, the values of the chemotactic
sensitivity coefficient $\chi$ and the  chemicals diffusion constant
$k$ play significant roles to determine the behavior of solutions.
Obviously, either small $\chi>0$ or large $k>0$ benefits the global
boundedness of solutions.

Recall the known results in the field with $k=1$. At first consider
the parabolic-elliptic analogue of (\ref{P}), namely, the second
parabolic equation in (\ref{P}) is replaced by the elliptic equation
$0=\Delta v-v+u$. It was known that all radial classical solutions
are global-in-time if either $n\geq 3$ with $\chi<\frac{2}{n-2}$, or
$n=2$ with $\chi>0$ arbitrary  \cite{TT}. When
$0<\phi(v)<\frac{\chi}{v^l}$ with $l\geq1$, $\chi>0$, there is a
unique and globally bounded classical solution if $\chi<\frac{2}{n}$
$(l=1)$ or
$\chi<\frac{2}{n}\cdot\frac{l^l}{(l-1)^{(l-1)}}\gamma^{(l-1)}$
$(l>1)$, with $\gamma>0$ depending on $\Omega$ and $u_0$ \cite{KMT}.
Next consider the parabolic-parabolic case. All solutions of
(\ref{P}) are global in time when either $n=1$ \cite{KA}, or $n=2$
and $\chi<\frac{5}{2}$ under the radial assumption, while $\chi<1$
under the non-radial assumption \cite{P,TTK}. For $n\ge 2$,
(\ref{P}) possesses global classical solutions if
$0<\chi<\sqrt\frac{2}{n}$, and moreover,
$\chi<\sqrt{\frac{n+2}{3n-4}}$ ensures the global existence of weak
solutions \cite{M}. Certain radial weak solutions have been
constructed when $\chi<\sqrt{\frac{n}{n-2}}$  \cite{CN}. Following,
the global boundedness of solutions has been obtained for
$\chi<\sqrt{\frac{2}{n}}$  \cite{K}. Refer to \cite{KT,J,Z} for
the results on more general chemotaxis models.

Recently, under somewhat complicated conditions, Wang \cite{W}
established classical global solutions to the problem, a similar
model to \eqref{P},
\begin{eqnarray*} \left\{
 \begin{array}{llll}
  \displaystyle u_t=\nabla\cdot(\nabla u-\frac{\chi u}{v+c}\nabla v),& x\in\Omega,~~t\in(0,T),\\[4pt]
  \displaystyle v_t=k\Delta v-\alpha v+\beta u,& x\in\Omega,~~t\in(0,T),\\[4pt]
  \displaystyle \frac{\partial u}{\partial {\nu}}=\frac{\partial v}{\partial {\nu}}=0 ,& x\in\partial\Omega,~~t\in(0,T),\\[4pt]
  \displaystyle (u(x,0),v(x,0))=(u_0(x),v_0(x)),  &x\in\Omega,
 \end{array}\right.
\end{eqnarray*}
with $\chi,c,k,\alpha,\beta>0$. In the present paper, motivated by
Winkler \cite{M} and Fujie \cite{K}, we will prove the global
existence-boundedness of classical solutions to \eqref{P}, with
simplified conditions. That is the following theorem.
\begin{theorem}\label{the2}
Let $n\geq 2$, $u_0\in C^0(\overline\Omega)$, $v_0\in
W^{1,q}(\Omega)$ $(q>n)$ with $u_0\ge 0$, $v_0>0$ on $\overline\Omega$.
 Then, for any $k>0$, there exists a global classical solution to
\eqref{P}, provided
$\chi\in\big(0,-\frac{k-1}{2}+\frac{1}{2}\sqrt{(k-1)^2+\frac{8k}{n}}\big)$.
Moreover, the solution is globally bounded under $n\leq 8$.
\end{theorem}
\begin{remark}\label{R1} {\rm
Theorem \ref{the2} shows in what way the size of $k>0$ (the
diffusion strength of the chemicals $v$) effects the behavior of
solutions to \eqref{P}. It is interesting to observe that when $n=2$
the global existence-boundedness of solutions is independent of the
size of $k>0$, since
$-\frac{k-1}{2}+\frac{1}{2}\sqrt{(k-1)^2+\frac{8k}{n}}\equiv 1$ with
$n=2$. Quite differently, when $n\ge 3$ the contribution of $k>0$ is
significant that the range of $\chi$ for global
existence-boundedness of solutions will be enlarged (shrunk) as
$k>0$ is degreasing (increasing). The arbitrariness of $k>0$ yields
the ``maximal" range with $\chi\in(0,1)$ or the ``minimal" range
with $\chi\in (0,\frac{2}{n})$. That is to say for any
$\chi\in(0,1)$ (close to 1), there is $k>0$ (small) such that the
classical solution of \eqref{P} is globally bounded. On the other
hand, for any $k>0$ (large), there is $\chi\in (0,\frac{2}{n})$ to
ensure the global boundedness. Finally, it is pointed out that if
$k=1$, the required range of $\chi$ in Theorem \ref{the2} becomes
$0<\chi<\sqrt{\frac{2}{n}}$, which coincides with those in
\cite{K,M}}.
\end{remark}

\begin{remark} \label{R2}{\rm
Now compare our results for the parabolic-parabolic chemotaxis model
\eqref{P} with those for the corresponding parabiotic-elliptic
model, which can be thought as a special case of \eqref{P} with the
diffusion constant of the chemicals $v$ sufficiently large
\cite{KMT}. Just as mentioned in Remark \ref{R1}, letting $k>0$ be
arbitrarily large results in the ``minimal" permitted range with
$0<\chi<\frac{2}{n}$.
 This does agree those obtained for the
parabolic-elliptic model in \cite{TT}}.
\end{remark}

\section{Preliminaries}

In this section we introduce the local existence of classical
solutions to \eqref{P} with required estimates involving $\chi$ and
$k$, as well as some technical lemmas for the global boundedness as
preliminaries.

\begin{lemma}\label{lm2.1}
 Let $n\geq 2$, $u_0\in C^0(\overline\Omega)$, $v_0\in
W^{1,q}(\Omega)$ $(q>n)$ with $u_0\ge 0$, $v_0>0$ on $\overline\Omega$.
Then, for any $k,\chi>0$, there exists $T_{\max}\in (0,\infty]$, such that (\ref{P}) has
     a unique nonnegative solution $u\in C^0([0,T_{\max});C^0(\Omega))\times C^{2,1}
     (\overline{\Omega}\times(0,T_{\max}))$ and $v\in C^0([0,T_{\max});C^0(\Omega))\times C^{2,1}
     (\overline{\Omega}\times(0,T_{\max}))\times L_{\rm loc}^{\infty}\big([0,T);~W^{1,q}(\Omega)\big)$, where
      either $ T_{\max}=\infty$, or $T_{\max}<\infty $ with
  $\lim_{t\rightarrow T_{\max}}\| u(\cdot,t)\|_{L^\infty(\Omega)}+\| v(\cdot,t)\|_{W^{1,q}(\Omega)}=\infty.$
\end{lemma}
{\bf Proof.}\ For $k>0$, it is known from Lemma 2.2 in \cite{K} that
there is $\eta>0$, such that $\inf_{x\in\Omega}v(x,t)\ge \eta>0$ for
all $t>0$. Consequently,  the local existence lemma can be obtained
by the classical parabolic theory, refer to Theorem 3.1 in
\cite{MH}.\qquad$\Box$\medskip

 The following lemma is crucial to
establish the global existence-boundedness conclusions of the paper.
Denote $r_{\pm}(p)=(p-1)[\frac{p\chi(1-k)+2k}{p(k-1)^2+4k}
\pm\frac{2\sqrt{k^2-p\chi k(k-1)-p\chi^2k}}{p(k-1)^2+4k}]$.
Throughout the paper, for simplicity, denote $T=T_{\max}$.

\begin{lemma}\label{lemma24}
 Let $(u,v)$ solve \eqref{P} with $k>0$ and
  $\chi\in\big(0,-\frac{k-1}{2}+\frac{1}{2}\sqrt{(k-1)^2+\frac{8k}{n}}\big)$.
  If
\begin{align}\label{0}
\|v(\cdot,t)\|_{L^{p-r}(\Omega)}\leq c, ~~t\in(0,T)
\end{align}
with $p<\frac{k}{[\chi^2-\chi(1-k)]_+}$,  $r\in
(r_{-}(p),r_{+}(p))$, and $c>0$,  then
\begin{align}\label{33}\int_{\Omega}u^{p}v^{-r}dx\leq
\tilde{c},~~t\in(0,T) \end{align} with some $\tilde{c}>0$.
\end{lemma}
{\bf Proof.}\
  It is known via a simple computation with (\ref{P}) that
\begin{align*}
  \frac{d}{dt} \int_\Omega u^{p}v^{-r} dx &=p\int_\Omega u^{p-1}v^{-r}[\Delta u-\chi\nabla(\frac{u}{v}\nabla v)] dx-r\int_\Omega u^p v^{-r-1}(k\Delta v-v+u) dx\\
   &=-p\int_\Omega \nabla(u^{p-1}v^{-r})\cdot(\nabla u-\chi\frac{u}{v}\nabla v) dx+rk\int_\Omega\nabla (u^pv^{-r-1})\cdot\nabla v dx\\
   &~~~~~ +r\int_\Omega u^pv^{-r} dx -r\int_\Omega u^{p+1}v^{-r-1} dx \\
    &=-p(p-1)\int_\Omega u^{p-2}v^{-r}|\nabla u|^2 dx +[pr+prk+p(p-1)\chi]\int_\Omega u^{p-1}v^{-r-1}\nabla u\cdot\nabla  v dx\\
    &~~~~-[r(r+1)k+pr\chi]\int_\Omega u^pv^{-r-2}|\nabla v|^2 dx
    +r\int_\Omega u^pv^{-r} dx -r\int_\Omega u^{p+1}v^{-r-1}dx\\
 &\leq\int_\Omega \Big[{\frac{p[(p-1)\chi+r+rk]^2}{4(p-1)}
-pr\chi-r(r+1)k}\Big]u^pv^{-r-2}|\nabla v|^2 dx\\
 &~~~~~+r\int_\Omega u^pv^{-r} dx -r\int_\Omega u^{p+1}v^{-r-1} dx
 \end{align*}
by Young's inequality. Denote
\begin{align*}
 f(r;p,\chi,k)=\frac{p[(p-1)\chi+r+rk]^2}{4(p-1)}
 -pr\chi-r(r+1)k,
 \end{align*}
and rewrite as the quadric expression in $r$
\begin{align*}
4(p-1)f(r;p,\chi,k)=[p(k-1)^2+4k]r^2+[2p(p-1)\chi(k-1)-4(p-1)k]r+p(p-1)^2
\chi^2.
\end{align*}
We know
\begin{align*}
\Delta_r&=4(p-1)^2[p\chi(k-1)-2k]^2-4(p-1)^2p\chi^2[p(k-1)^2+4k]\\
&=16(p-1)^2[k^2-p\chi k(k-1)-p\chi^2k]>0
\end{align*}
whenever $ p<\frac{k}{[\chi^2+\chi(k-1)]_+}$. Consequently,
$f(r;p,\chi,k)< 0$ for any $r\in \big(r_{-}(p),r_{+}(p)\big)$. This
yields
 \begin{align}\label{ineq23}
\frac{d}{dt}\int_\Omega u^{p}v^{-r}dx\leq
r\int_\Omega u^pv^{-r}dx-r\int_\Omega
u^{p+1}v^{-r-1}dx,  ~~t\in(0,T).
\end{align}
Due to $\int_\Omega u^pv^{-r}\leq\big(\int_\Omega
u^{p+1}v^{-r-1}\big)^\frac{p}{p+1}\big(\int_\Omega
v^{p-r}\big)^\frac{1}{p+1}$, we obtain \eqref{33} from
\eqref{ineq23} and  (\ref{0}).\qquad $\Box$\medskip

\begin{lemma} \label{l22}
Let $(u,v)$ satisfy the second equation of \eqref{P} with $k>0$,
$1\leq q\leq p\leq\infty$, $\frac{n}{2}(\frac{1}{q}-\frac{1}{p})<1$.
Then there exists $C>0$, such that
\begin{align}\label{ineq21}
\|v(\cdot,t)\|_{L^p(\Omega)}\leq
C(1+\sup_{s\in(0,t)}\|u(\cdot,s)\|_{L^q(\Omega)}),~~ t\in(0,T)
\end{align}
\end{lemma}
{\bf Proof.} Noticing $v(\cdot,t)={\rm
e}^{t(k\Delta-1)}v_0+\int_0^t{\rm e}^{(t-s)(k\Delta-1)}u(\cdot,s)ds$
for $t>0$, by the standard smoothing estimates for the heat
semigroup under homogeneous Neumann boundary conditions \cite{MW},
we can obtain for $q\leq p$ that
\begin{align*}
&\|v(\cdot,t)\|_{L^p(\Omega)}\le \|{\rm
e}^{t(k\Delta-1)}v_0\|_{L^p(\Omega)}+\int_0^t\|{\rm e}^{(t-s)(k\Delta-1)}u(\cdot,s)\|_{L^p(\Omega)}ds\\
&\le
C_1\|v_0(x)\|_{L^{\infty}(\Omega)}+C_2\sup_{s\in(0,t)}\|u(\cdot,s)\|_{L^q(\Omega)}\int_0^t(1+[k(t-s)])
^{-\frac{n}{2}(\frac{1}{q}-\frac{1}{p})}{\rm e}^{-(\lambda_1+\frac{1}{k})[k(t-s)]}ds\\
&\le
C_1\|v_0(x)\|_{L^{\infty}(\Omega)}+\frac{C_2}{k}\sup_{s\in(0,t)}\|u(\cdot,s)\|_{L^q(\Omega)}\int_0^{\infty}(1+\alpha)
^{-\frac{n}{2}(\frac{1}{q}-\frac{1}{p})}{\rm e}^{-\lambda_1\alpha}
d\alpha,~~t\in(0,T),
\end{align*}
where $C_1, C_2>0$, $\lambda_1$ is the first nonzero eigenvalue of
$-\Delta$ under the Neumann boundary condition. This proves
\eqref{l22}.\qquad$\Box$
\medskip

Instead of the choice of $r=\frac{p-1}{2}$ in \cite{K}, we deal with
the more complicated form
$r=\small(p-1\small)\frac{[p\chi(1-k)+2k]}{p(1-k)^2+4k}$ to describe
the effect of the chemicals diffusion rate $k$, and denote
$h(p)\triangleq h(p;\chi,k)=\frac{p\chi(1-k)+2k}{p(1-k)^2+4k}$. For
such $h(p)$, we have the following lemma.
 \begin{lemma}\label{lemma 2.4}
Let $k>0$,
$\chi\in\big(0,-\frac{k-1}{2}+\frac{1}{2}\sqrt{(k-1)^2+\frac{8k}{n}}\big)$,
and $p\in\big(1,\frac{k}{[\chi^2+\chi(k-1)]_+}\big)$. Then
$h(p)\in(0,1)$.
\end{lemma}
{\bf Proof.} We have
$h'(p)=\frac{2k(1-k)[2\chi-(1-k)]}{[p(1-k)^2+4k]^2}$.\smallskip

If $k=1$, then $h(p)\equiv \frac{1}{2}$.

If $k>1$, then $h'(p)<0$, and so\\
$0<\frac{\chi}{2\chi+k-1}=h(\frac{k}{\chi^2-\chi(1-k)}-0)
=h(\frac{k}{[\chi^2-\chi(1-k)]_+}-0)<h(p)<h(1+0)
=\frac{2k+\chi(1-k)}{(1-k)^2+4k}<\frac{2k}{(1+k)^2}<1$.\medskip

Now suppose $0<k<1$.

If $0<\chi<\frac{1-k}{2}$, then $h'(p)<0$, and so\\
$0<\frac{\chi}{1-k}=\mathop{\lim}\limits_{p\rightarrow\infty}
\frac{p\chi(1-k)+2k}{p(1-k)^2+4k}<h(p)<h(1+0)
=\frac{2k+\chi(1-k)}{(1-k)^2+4k} <\frac{1}{2}$;\medskip

If $\chi=\frac{1-k}{2}$, then $h(p)\equiv \frac{1}{2}$;\medskip

If $\frac{1-k}{2}<\chi\leq 1-k$, then $h'(p)>0$, and so \\
$\frac{1}{2}<\frac{2k+\chi(1-k)}{(1+k)^2}=h(1+0)<h(p)
<\mathop{\lim}\limits_{p\rightarrow\infty}\frac{p\chi(1-k)+2k}
{p(1-k)^2+4k}=\frac{\chi}{1-k}\leq 1$;\medskip

If $1-k<\chi<\frac{-(k-1)+\sqrt{(k-1)^2+\frac{8k}{n}}}{2}$, then
$h'(p)>0$, and so\\
$0<\frac{1-k^2}{(1+k)^2}<\frac{2k+\chi(1-k)}{(1+k)^2}=h(1+0)<h(p)
<h(\frac{k}{\chi^2-\chi(1-k)}-0)
=\frac{\chi}{2\chi+k-1}<1$.\medskip

The proof is complete.\qquad$\Box$\medskip

Denote $c_0=\inf_{p\in(1,\frac{n}{2}]}h(p)$ and
$c^0=\sup_{p\in(1,\frac{n}{2}]}h(p)$ with $n\ge 3$. By Lemma
\ref{lemma 2.4} and its proof, $c_0,c^0\in (0,1)$. The following
lemma with $c_0$ and $c^0$  will play an important role for
estimating the bound of $u$ by the involved iteration in the next
section.
\begin{lemma}\label{cor3.1}
Let $k>0$,
$\chi\in\big(0,-\frac{k-1}{2}+\frac{1}{2}\sqrt{(k-1)^2+\frac{8k}{n}}\big)$,
then
$f(x)\triangleq\frac{n[(1-c^0)x+c^0-c_0]+2c_0x}{(n-2x)(1-c_0)}-x>0$
for $x\in(1,\frac{n}{2}]$, provided $n\le 8$.
\end{lemma}
{\bf Proof.} It suffices to show
$g(x)\triangleq2(1-c_0)x^2+[2c_0-n(c^0-c_0]x+n(c^0-c_0)>0$ in
$(1,\frac{n}{2}]$.

The case of $c_0=c^0\in(0,1)$ is trivial.

Now suppose $0<c_0<c^0<1$. We have
$\Delta_g=[2c_0-n(c^0-c_0)]^2-8n(1-c_0)(c^0-c_0)
=[n(c^0-c_0)]^2+(4c_0-8)n(c^0-c_0)+4c_0^2<0$, and hence $g(x)>0$,
whenever $4-2c_0-4\sqrt{1-c_0}< n(c^0-c_0)<4-2c_0+4\sqrt{1-c_0}$.

If $n(c^0-c_0)\leq 4-2c_0-4\sqrt{1-c_0}$, then $\Delta_g>0$. Due
$n(c^0-c_0)-2c_0<4(1-c_0)-4\sqrt{1-c_0}<0$, we know that both the
two roots of $g(x)$ must be negative. With $g(1+0)= 2$, we obtain
$g(x)>0$.

If $n(c^0-c_0)\geq 4-2c_0+4\sqrt{1-c_0}$, then $\Delta_g>0$, and the
two roots of $g(x)$ satisfy $x_2\ge x_1>0$. Together with $g(1+0)=2$
and $g(\frac{n}{2})=\frac{n^2}{2}(1-c^0)+nc_0>0$, the positivity of
$g(x)$ for $x\in (1,\frac{n}{2}]$ requires that the minimal point of
$g(x)$ satisfies $\frac{n(c^0-c_0)-2c_0}{4(1-c_0)}<1$, i.e.,
$n(c^0-c_0)+2c_0<4$, by Vieta's formulas. This contradicts the case
$n(c^0-c_0)\geq 4-2c_0+4\sqrt{1-c_0}$. So, the case itself should be
excluded to ensure $g(x)>0$ in $(1,\frac{n}{2}]$.   Rewrite the case
as $c^0\ge \frac{(n-2)c_0+4+4\sqrt{1-c_0}}{n}\triangleq \alpha(c_0)$
with $c_0\in(0,1)$. We have $c^0\ge
\alpha(c_0)>\min\{\alpha(0),\alpha(1)\}
=\min\{\frac{n+2}{n},\frac{8}{n}\}$. Thus, we would get a
contradiction that  $c^0>1$, whenever $n\le 8$. \qquad$\Box$

\section{Proof of main result}

We deal with the proof of the main result of the paper in this
section. \\[6pt]
 {\bf Proof of Theorem \ref{the2}.}

We at first show that the local solutions ensured by Lemma
\ref{lm2.1} should be global.

 Assume $p>q$. By the H\"{o}lder
inequality, we have
\begin{align}\label{eq31}
 \int_\Omega u^qdx&=\int_\Omega (u^pv^{-r})^{\frac{q}{p}}v^{\frac{rq}{p}}dx
 \leq\big(\int_\Omega u^pv^{-r}dx\big)^{\frac{q}{p}}\big(\int_\Omega v^{\frac{rq}{p-q}}dx\big)
 ^{\frac{p-q}{p}}, ~~0<t<T.
 \end{align}

Let  $ p<\frac{k}{[\chi^2+\chi(k-1)]_+}$. We know from
\eqref{ineq23} in the proof of Lemma \ref{lemma24} with $r\in
(r_{-}(p),r_{+}(p))$ that $$ \frac{d}{dt}\int_\Omega
u^{p}v^{-r}dx\leq r\int_\Omega u^pv^{-r}dx, ~~t\in(0,T),
$$ and hence
\begin{align}\label{00}\int_{\Omega}u^{p}v^{-r}dx\leq C,~~t\in(0,T)
\end{align} with $C=C(t)>0$.  Notice $\chi\in\big(0,-\frac{k-1}{2}+\frac{1}{2}\sqrt{(k-1)^2+\frac{8k}{n}}\big)$ with $
p<\frac{k}{[\chi^2+\chi(k-1)]_+}$ admits $p>\frac{n}{2}$.

Take $q\in(\frac{n}{2},p)\subset
\big(1,\min\{p,\frac{n(p-r)}{[n-2r]_+}\}\big)$. Then
$\frac{n}{2}(\frac{1}{q}-\frac{p-q}{rq})<1$. By Lemma
 \ref{l22}, we have with $C_1>0$ that
\begin{align}\label{ineq33}
\| v\|_{\frac{rq}{p-q}}\leq C_1(1+\mathop{\sup}\limits_{s\in(0,t)}\|
u\|_{q}), ~~0<t<T.
\end{align}
Combining \eqref{eq31}--\eqref{ineq33} yields
\begin{align*}
\int_\Omega u^qdx\leq C_2\big(1+(\sup_{ s\in(0,t)}\int_\Omega
u^qdx)^{\frac{r}{p}}\big), ~~t\in(0,T),
\end{align*}
and hence $\int_\Omega u^qdx\leq C_3$ for $t\in(0,T)$, where
$C_2=C_2(t)=\tilde{C}_2 {\rm e}^{rt}$ with $\tilde{C}_2>0$, $
C_3=C_3(t)>0$, and $\frac{r}{p}<\frac{p-1}{p}<1$ by Lemma \ref{lemma
2.4}.
 Now we can follow the proof of  Lemma 3.4 in \cite {M} to obtain the global existence of solutions for \eqref{P}.
 \medskip

 Next, we prove the solutions established above are also globally
bounded if $n\leq 8$. We should verify that $\int_\Omega u^q dx< C$
with some $q>\frac{n}{2}$ and $C>0$ \cite{M}.

The case of $n=2$ is simple. Take
\begin{equation*}
 \left\{ \begin{aligned}
 &p\in\big(1,\frac{k}{[\chi^2+\chi(k-1)]_+}\big),\\
 &r=(p-1)h(p).
 \end{aligned} \right.
 \end{equation*}
Then $1-\frac{1}{p-r}<1$. By Lemma \ref{l22} with the
$L^1$-conservation of $u$, we have
 $$\|v\|_{L^{p-r}}\le C_4,~~t\in(0,T)$$ with time-independent $C_4>0$.
  By Lemma \ref{lemma24},
 \begin{align}\label{111}\int_{\Omega}u^{p}v^{-r}dx\leq
C_5,~~t\in(0,T) \end{align} with some $C_5>0$.  We deduce
 from \eqref{eq31}, \eqref{ineq33} and \eqref{111} that there exists $q>1$, such that $\| u\|_{q}\leq C_6$ for
 $t\in (0,T)$ with $C_6=C_6(k)>0$.

Now consider the case of $3\le n\leq 8$.
 We should fix the assumption \eqref{0} with
some $p>\frac{n}{2}$ and the constant independent of $t$ there. We
will do it via an iteration procedure based of Lemma \ref{cor3.1}.

$1^\circ$  Take
\begin{equation}\label{11}
 \left\{ \begin{aligned}
 &p_0\in\Big(1,\min\big\{\frac{k}{[\chi^2+\chi(k-1)]_+},
 \frac{n(1-c_0)+2c_0}{(n-2)(1-c_0)}\big\}\Big),\\
 &r_0=(p_0-1)h(p_0).
 \end{aligned} \right.
 \end{equation}
 Then $p_0-r_0\le p_0-( p_0-1)c_0=(1-c_0)p_0+c_0<\frac{n}{n-2}$,
 i.e., $\frac{n}{2}(1-\frac{1}{p_0-r_0})<1$. By Lemma \ref{l22} with the $L^1$-conservation of $u$, we have
 $$\|v\|_{L^{p_0-r_0}}\le C_7,~~t\in(0,T),$$and hence
 \begin{align}\label{333}\int_{\Omega}u^{p_0}v^{-r_0}dx\leq
C_8,~~t\in(0,T) \end{align} by Lemma \ref{lemma24}, with some
$C_7,C_8>0$.  We know
 from \eqref{eq31}, \eqref{ineq33} and \eqref{333} that $\| u\|_{q_0}\leq C_9$ for
 $t\in (0,T)$ with $C_9>0$. If
$p_0>\frac{n}{2}$, take $q_0\in(\frac{n}{2},p_0)\subset
 \big(1,\min\{p_0,\frac{n(p_0-r_0)}{[n-2r_0]_+}\}\big)$.\medskip

$2^\circ$     Assume $p_0\leq\frac{n}{2}$. Take
\begin{equation*}
\left\{ \begin{aligned}
&p_1\in\Big(p_0,\min\big\{\frac{k}{[\chi^2+\chi(k-1)]_+},
\frac{n[(1-c^0)p_0+c^0-c_0]+2c_0p_0}{(n-2p_0)(1-c_0)}\big\}\Big),\\
&r_1=(p_1-1)h(p_1),
 \end{aligned} \right.
 \end{equation*}
 where the well-definedness of the interval for $p_1$ is ensured  by Lemma
 \ref{cor3.1}.
A simple calculation shows \begin{align*}
p_1-r_1\leq(1-c_0)p_1+c_0<\frac{n[(1-c^0)p_0+c^0]}{n-2p_0}\leq
\frac{n(p_0-r_0)}{n-2p_0}=\frac{n\frac{n(p_0-r_0)}{n-2r_0}}
{n-\frac{2n(p_0-r_0)}{n-2r_0}},
\end{align*}
i.e.,$\frac{n}{2}(\frac{1}{\frac{n(p_0-r_0)}{n-2r_0}}-\frac{1}{p_1-r_1})<1$.
By Lemma \ref{l22} with $p_0\le\frac{n}{2}$, there is
$q_0\in(1,\frac{n(p_0-r_0)}{n-2r_0})$ such that $$\|
v\|_{p_1-r_1}<C_{10}, ~~ t\in(0,T),$$ and thus
\begin{align}\label{444}
\int_{\Omega}u^{p_1}v^{-r_1}dx\leq C_{11},~~t\in(0,T) \end{align} by
Lemma \ref{lemma24}, with some $C_{10},C_{11}>0$. It follows from
\eqref{eq31}, \eqref{ineq33} and \eqref{444} that $\|u\|_{q_1}\le
C_{12}$ for $t\in(0,T)$ with $C_{12}>0$. If $p_1>\frac{n}{2}$, take
$q_1\in(\frac{n}{2},p_1)\subset
 \big(1,\min\{p_1,\frac{n(p_1-r_1)}{[n-2r_1]_+}\}\big)$.\medskip

$3^\circ$  Assume $p_{l-1}\leq\frac{n}{2}$, $l=1,2,3,\dots$. Take
\begin{equation}\label{12}
\left\{ \begin{aligned}
&p_l\in\Big(p_{l-1},\min\big\{\frac{k}{[\chi^2+\chi(k-1)]_+},
\frac{n[(1-c^0)p_{l-1}+c^0-c_0]+2c_0p_{l-1}}{(n-2p_{l-1})
(1-c_0)}\big\}\Big),\\
&r_l=(p_{l-1}-1)h(p_{l-1}),
 \end{aligned} \right.
 \end{equation}
  where the interval for $p_l$ is well defined due to Lemma
 \ref{cor3.1}.
Repeat the procedure in $2^\circ$, we deduce with
$q_l\in\big(1,\min\{p_l,\frac{n(p_l-r_l)}{[n-2r_l]_+}\}\big)$ that
 $\| u\|_{q_l}\leq C_{13}$, $0<t<T$, with some $C_{13}>0$.
Noticing $\frac{n[(1-c^0)p_{l-1}+c^0-c_0]+2c_0p_{l-1}}{(n-2p_{l-1})
(1-c_0)}\rightarrow \infty$ as $l\rightarrow \infty$ by Lemma
 \ref{cor3.1}, we can realize $p_l>\frac{n}{2}$ after finite steps.
 Let $q=q_l\in(\frac{n}{2},p_l)$ to get $\int_\Omega u^qdx \le C_{14}$
  for all $t\in (0,T)$ with $C_{14}>0$. It is mentioned that for any fixed $k>0$, the
  constants $C_i$, $i=5,\dots,14$, are all independent of
  $t\in(0,T)$ here.

Based on the above estimate $\|u(\cdot,t)\|_{L^q(\Omega)}$ with
$q>\frac{n}{2}$, uniform for $t\in (0,T)$, we conclude the global
boundedness of $u$ by coping the related arguments in \cite{M} for
that to the case of $k=1$.\qquad $\Box$
\begin{remark}\label{R3} {\rm
 Notice $
p_0<\frac{k}{[\chi^2+\chi(k-1)]_+}$ in \eqref{11} with
$\chi\in\big(0,-\frac{k-1}{2}+\frac{1}{2}\sqrt{(k-1)^2+\frac{8k}{n}}\big)$
admits $p_0>\frac{n}{2}$. Moreover, a simple computation shows
$p_0<\frac{n(1-c_0)+2c_0}{(n-2)(1-c_0)}$ in \eqref{11}  ensures
$\frac{n(1-c_0)+2c_0}{(n-2)(1-c_0)}>\frac{n}{2}$ whenever $n=3,4$.
Therefore, the Steps $2^\circ$ and $3^\circ$ in the proof of Theorem
\ref{the2} are unnecessary for $n=3,4$ there. In addition, it should
be pointed out that $c^0=c_0=\frac{1}{2}$ if $k=1$, and thus the
requirement $n\le 8$ itself can be removed away for the global
boundedness of solutions with $k=1$. }\end{remark}


\end{document}